\documentclass[reqno,10pt]{article}

\usepackage{amsmath}
\usepackage{amssymb}
\usepackage{amsthm}
\usepackage[UKenglish]{babel}
\usepackage{bbm}               % for tfn
\usepackage{ifthen}
\usepackage{sectsty}
\usepackage{url}

\sectionfont{\large}
\subsectionfont{\normalsize}

\mathchardef\ordinarycolon\mathcode`\:
\mathcode`\:=\string"8000
\def\vcentcolon{\mathrel{\mathop\ordinarycolon}}
\begingroup \catcode`\:=\active
  \lowercase{\endgroup
  \let :\vcentcolon
  }

\newtheoremstyle{mythm}{3pt}{3pt}{\itshape}{}{\bfseries}{}{ }%
{\thmname{#1}~\thmnumber{#2}.\thmnote{~(#3)}}

\theoremstyle{mythm}
\newtheorem{theorem}{Theorem}[section]
\newtheorem{lemma}{Lemma}[section]
\newtheorem{proposition}{Proposition}[section]

\theoremstyle{definition}
\newtheorem{definition}{Definition}[section]
\newtheorem{notation}{Notation}[section]

\newtheorem{remark}{Remark}[section]

\let\origthebibliography\thebibliography
\def\otherbibliography{\renewcommand{\section}[2]{}\origthebibliography}
\renewcommand{\thebibliography}[1]{\otherbibliography{#1}%
\setlength{\itemsep}{0.5ex}}

\newcommand{\ag}{\mathcal{D}}
\newcommand{\bebe}{\Gamma}
\newcommand{\bbebe}{{\widetilde{\bebe}}}
\newcommand{\di}{\mbox{$\Sigma$}}
\newcommand{\eehm}{\widetilde{\ehm}}
\newcommand{\eevecs}{\widetilde{\evecs}}
\newcommand{\ehm}{s}
\newcommand{\elltwo}{\mathsf{K}}
\newcommand{\expn}{\mathbb{E}}
\newcommand{\evec}[1]{\evecc(#1)}
\newcommand{\evecc}{\varepsilon}
\newcommand{\evecs}{\mathcal{E}}
\newcommand{\fflip}{\widetilde{\flip}}
\newcommand{\ffock}{{\widetilde{\fock}}}
\newcommand{\flip}{\sigma}
\newcommand{\fock}{\mathcal{F}}
\newcommand{\hilb}{\mathsf{H}}
\newcommand{\id}{I}
\newcommand{\indf}[1]{\indff_{#1}}
\newcommand{\indff}{1}
\newcommand{\ini}{\mathsf{h}}
\newcommand{\Ito}{\mathcal{I}}
\newcommand{\mul}{\mathsf{k}}
\newcommand{\mmul}{{\widehat{\mul}}}
\newcommand{\nind}[1]{\widetilde{\indff}_{#1}}
\newcommand{\Ptn}{\mathrm{T}}
\newcommand{\ptn}{\tau}
\newcommand{\Vac}{\Omega}
\newcommand{\vac}{\omega}
\newcommand{\vint}{\Lambda_\Omega}

\renewcommand{\ge}{\geqslant}
\renewcommand{\le}{\leqslant}

\newcommand{\bp}{\mathbf{p}}
\newcommand{\bq}{\mathbf{q}}
\newcommand{\bt}{\mathbf{t}}

\newcommand{\C}{\mathbb{C}}
\newcommand{\R}{\mathbb{R}}
\newcommand{\Z}{\mathbb{Z}}

\newcommand{\algten}{\odot}
\newcommand{\bigalgten}{\bigodot}
\newcommand{\bop}[2]%
{\ifthenelse{\equal{#2}{}}{\bopp(#1)}{\bopp(#1;#2)}}
\newcommand{\bopp}{\mathcal{B}}
\newcommand{\lin}{\mathop{\mathrm{lin}}}
\newcommand{\lop}[2]%
{\ifthenelse{\equal{#2}{}}{\mathcal{L}(#1)}{\lopp(#1;#2)}}
\newcommand{\lopp}{\mathcal{L}}
\newcommand{\mvec}[2]%
{\left(\begin{smallmatrix}#1\\[0.5ex]#2\end{smallmatrix}\right)}
\newcommand{\nnabla}{\widehat{\nabla}}
\newcommand{\rd}{\mathrm{d}}
\newcommand{\std}{\,\rd}
\newcommand{\tfn}[1]{\mathbbm{1}_{#1}}

\newenvironment{myquotation}{\small\quad}{\par\vspace{1ex}}

\newenvironment{mylist}%
{\begin{list}{}%
{\leftmargin 4.75em\labelwidth 3.5em\rightmargin 1.25em%
\topsep 0.75ex\itemsep 0.25ex}}%
{\end{list}}

\begin{document}

\begin{center}
{\LARGE
Approximation \emph{via} toy Fock space\\[0.5ex]
-- the vacuum-adapted viewpoint}\\[3ex]
{\large Alexander C. R. Belton}\\[2ex]
{\small School of Mathematical Sciences,\\
University College, Cork, Ireland.\\[1ex]
\textsf{a.belton@ucc.ie}}\\[2ex]
\end{center}

\begin{center}
\textit{Dedicated to Slava Belavkin on the occasion of his 60th birthday}
\end{center}

\begin{abstract}
After a review of how Boson Fock space (of arbitrary multiplicity) may be
approximated by a countable Hilbert-space tensor product (known as toy Fock
space) it is shown that vacuum-adapted multiple quantum Wiener integrals of
bounded operators may be expressed as limits of sums of operators defined on
this toy space, with strong convergence on the exponential domain. The
vacuum-adapted quantum It\^o product formula is derived with the aid of this
approximation and a brief pointer is given towards the unbounded case.
\end{abstract}

\section{Introduction}

The idea of using discrete approximations in quantum stochastic calculus goes
back at least as far as Meyer's notes \cite{Mey86}, where he gave credit to
Jean-Lin Journ\'e. Around the same time, articles by Parthasarathy
\cite{Par88} and Lindsay and Parthasarathy \cite{LiP88} showed that
certain quantum flows (which are generalisations of classical diffusions) may
be approximated by so-called spin random walks, while Accardi and Bach
produced (in an unpublished preprint --- see \cite{AcB89,Mey89}) a
central-limit theorem which may be viewed as a result on toy-Fock-space
approximation. These ideas have recently been the subject of renewed interest.

Attal revisited and extended the Journ\'e--Meyer ideas in \cite{Att03}, giving
a heuristic derivation of the quantum It\^o product formula using the
approximation, and this was followed by further work of Attal and Pautrat
\cite{AtP06} and of Pautrat \cite{Pau05}. Their point of view may be
considered as physical, rather than probabilistic; in \cite{Gou04}, Gough
examined the physics of this set-up and explained its connexion with Holevo's
time-ordered exponentials.

Meanwhile, Sinha \cite{Sin06} revived the ideas of \cite{LiP88},
emphasising that, in many cases, sufficiently strong convergence holds to
enable one to deduce that the limit flows are $*$-homomorphic. Further work in
this direction has done by Sahu \cite{Sah05}, who moved away from the spin
approach of Lindsay--Parthasarathy--Sinha to adopt the same type of coupling
between system and noise as Attal--Pautrat; it is this direct (as opposed to
spin) coupling which is used below.

Many other people have worked with these concepts, including Bouten,
van~Handel and James \cite{BHJ06} (in quantum filtering), Brun
\cite{Bru02} and Gough and Sobolev \cite{GoS04} (who view the situation as
physicists), Franz and Skalski \cite{FrS07} (for constructing random walks
on quantum groups), K\"ummerer (\cite{Kum06} gives a detailed physical
interpretation of discrete models and is an excellent introduction to his
earlier work in this area) and Leitz-Martini \cite{Lei01} (who expressed
many of these approximation ideas rigorously using non-standard analysis; for
example, the discrete It\^o table of Attal \cite[Section~VII]{Att03} agrees
with the continuous-time version only in the limit, but in the non-standard
setting the anomalous terms are infinitesimal \cite[(2.2.8)]{Lei01}).

Here, a vacuum-adapted approach to approximation is adopted and, as might be
expected, a very straightforward theory results. After revising the embedding
of toy Fock space into Boson Fock space in Section~\ref{sec:space}, modified
versions of the vacuum-adapted Wiener integral are defined in
Section~\ref{sec:modint}. The natural `discrete integral' (which is, of
course, a sum) is examined in Section~\ref{sec:toyint} and is shown to be
given, up to an error term, by the modified integral previously
defined. Section~\ref{sec:mulint} extends this working to the case of multiple
integrals, Section~\ref{sec:itof} shows how the quantum It\^o product formula
arises naturally from the discrete approximation and Section~\ref{sec:future}
points the way to further developments involving unbounded
operators. Applications of these results will appear elsewhere \cite{Blt07}.

\subsection{Conventions and Notation}
All sesquilinear inner products are conjugate linear in the first variable. We
follow \cite{Lin05} for the most part, although the ordering of certain
objects is changed: for us, the initial space always appears first (the
`usual' convention, to quote Lindsay \cite[p.183]{Lin05}).

The vector space of linear operators between vector spaces $U$ and $V$ is
denoted by $\lop{U}{V}$, or $\lop{U}{}$ if $U$ and $V$ are equal; the identity
operator on a vector space $V$ is denoted by $\id_V$. The operator space of
bounded operators between Hilbert spaces $\hilb_1$ and $\hilb_2$ is denoted by
$\bop{\hilb_1}{\hilb_2}$, abbreviated to $\bop{\hilb_1}{}$ if $\hilb_1$ equals
$\hilb_2$. The tensor product of Hilbert spaces and bounded operators is
denoted by $\otimes$; the algebraic tensor product is denoted by
$\algten$. The restriction of a function $f$ to a set $A$ is denoted by
$f|_A$; the indicator function of $A$ is denoted by $\indf{A}$. If $P$ is a
proposition then the expression $\tfn{P}$ has the value $1$ if $P$ is true and
$0$ if $P$ is false.

\subsection*{Acknowledgements}
This is an expanded version of a talk given at the 27th Conference on Quantum
Probability and its Applications, which was held at the University of
Nottingham in July~2006; the hospitality of its organisers is acknowledged
with pleasure. Thanks are extended to Professor Martin Lindsay, for a remark
clarifying the proof of Theorem~\ref{thm:vacest} and for extensive comments on
a previous draft, which (hopefully) have led to a much improved presentation
of this material; thanks are also extended to Dr Stephen Wills, for many
helpful conversations on these ideas. The author is an Embark Postdoctoral
Fellow at University College Cork, funded by the Irish Research Council for
Science, Engineering and Technology.

\section{Toy and Boson Fock spaces}\label{sec:space}

\begin{myquotation}
Men more frequently require to be reminded than informed.

{\small\ -- Samuel Johnson, \textit{The Rambler}, No. 2 (1749--50).}
\end{myquotation}

\begin{notation}
Let $\mul$ be a complex Hilbert space (called the \emph{multiplicity space})
and let $\mmul:=\C\oplus\mul$ be its one-dimensional extension. Elements of
$\mmul$ will be written as column vectors, with the first entry a complex
number and the second a vector in $\mul$; if $x\in\mul$ then
$\widehat{x}:=\mvec{1}{x}$.
\end{notation}

\begin{definition}
\emph{Toy Fock space} is the countable tensor product
\begin{equation}
\bebe:=\bigotimes_{n=0}^\infty\mmul_{(n)}
\end{equation}
with respect to the stabilising sequence
$\bigl(\vac_{(n)}:=\mvec{1}{0})_{n=0}^\infty$, where $\mmul_{(n)}:=\mmul$
for each $n$; the subscript $(n)$ is used here and below to indicate the
relevant copy. (For information on infinite tensor products of Hilbert spaces,
see, for example \cite[Exercise~11.5.29]{KaR97}.)

For all $n\in\Z_+:=\{0,1,2,\ldots\}$, let
\begin{equation}
\bebe_{n)}:=\bigotimes_{m=0}^{n-1}\mmul_{(m)}\quad\mbox{and}\quad%
\bebe_{[n}:=\bigotimes_{m=n}^\infty\mmul_{(m)},
\end{equation}
where $\bebe_{0)}:=\C$. The identity $\bebe=\bebe_{n)}\otimes\bebe_{[n}$ is
the analogue of the continuous tensor-product structure of Boson Fock space.
\end{definition} 

\begin{notation}
For any interval $A\subseteq\R_+$, let $\fock_A$ denote Boson Fock space over
$L^2(A;\mul)$ and let $\fock:=\fock_{\R_+}$. For further brevity, let
$\elltwo=L^2(\R_+;\mul)$.

If $\ptn:=\{0=\ptn_0<\ptn_1<\cdots<\ptn_n<\cdots\}$ is a partition of $\R_+$
(so that $\ptn_n\to\infty$ as $n\to\infty$) then there exists an isometric
isomorphism
\begin{equation}
\Pi_\ptn:\fock\stackrel{\cong}{\longrightarrow}\fock_\ptn:=%
\bigotimes_{n=0}^\infty\fock_{[\ptn_n,\ptn_{n+1}[};\ %
\evec{f}\mapsto\bigotimes_{n=0}^\infty\evec{f|_{[\ptn_n,\ptn_{n+1}[}},
\end{equation}
where the tensor product is taken with respect to the stabilising sequence
$\bigl(\Vac_{[\ptn_n,\ptn_{n+1}[}:=%
\evec{0|_{[\ptn_n,\ptn_{n+1}[}}\bigr)_{n=0}^\infty$
and $\evec{g}$ denotes the exponential vector in $\fock_A$ corresponding to
the function $g\in L^2(A;\mul)$. The set of all such partitions of $\R_+$ is
denoted by $\Ptn$.
\end{notation}

\begin{definition}
For all $\ptn\in\Ptn$ and $n\in\Z_+$, define the natural isometry
\begin{equation}
j[\ptn]_n:\mmul\to\fock_{[\ptn_n,\ptn_{n+1}[};\ \mvec{\lambda}{x}\mapsto%
\lambda\Vac_{[\ptn_n,\ptn_{n+1}[}+x\nind{[\ptn_n,\ptn_n+1[},
\end{equation}
where $\nind{A}:=\indf{A}/\|\indf{A}\|_{L^2(\R_+)}$ is the normalised
indicator function of the interval $A\subseteq\R_+$, viewed as an element of
the one-particle subspace of $\fock_A$. These give an isometric embedding
\begin{equation}
J_\ptn:\bebe\to\fock_\ptn;\ \bigotimes_{n=0}^\infty\theta_n\mapsto%
\bigotimes_{n=0}^\infty j[\ptn]_n(\theta_n).
\end{equation}
Note that $Q_\ptn:=\Pi_\ptn^* J_\ptn J_\ptn^*\Pi_\ptn$ is an orthogonal
projection on $\fock$ and
\begin{equation}\label{eqn:embed}
J_\ptn^*\Pi_\ptn\evec{f}=\bigotimes_{n=0}^\infty\widehat{f_\ptn(n)}\qquad%
\forall\,f\in\elltwo,
\end{equation}
where
\begin{equation}
f_\ptn(n):=\frac{1}{\sqrt{\ptn_{n+1}-\ptn_n}}%
\int_{\ptn_n}^{\ptn_{n+1}}f(t)\std t\qquad\forall\,n\in\Z_+.
\end{equation}
\end{definition}

\begin{notation}
Let $\Ptn$ be the directed set of all partitions of $\R_+$, ordered by
inclusion; the expression `$f_\ptn\to f$ as $|\ptn|\to0$' means that the net
$(f_\ptn)_{\ptn\in\Ptn}$ converges to $f$. For all $\ptn\in\Ptn$, let
$P_\ptn\in\bop{\elltwo}{}$ be the orthogonal projection given by
\begin{equation}\label{eqn:cexp}
P_\ptn f:=\sum_{n=0}^\infty\frac{1}{\ptn_{n+1}-\ptn_n}%
\int_{\ptn_n}^{\ptn_{n+1}}f(t)\std t\,\indf{[\ptn_n,\ptn_{n+1}[}%
\qquad\forall\,f\in\elltwo.
\end{equation}
\end{notation}

\begin{lemma}\label{lem:cexp}
The projection $P_\ptn$ converges strongly to $\id_\elltwo$ as $|\ptn|\to0$.
\end{lemma}
\begin{proof}
If $f\in\elltwo$ is continuous and compactly supported, a uniform-continuity
argument may be used to show that $P_\ptn f\to f$ uniformly; the density of
such functions in $\elltwo$ completes the proof.
\end{proof}

\begin{theorem}
As $|\ptn|\to0$, the projection $Q_\ptn$ converges strongly to $\id_\fock$.
\end{theorem}
\begin{proof}
By (\ref{eqn:embed}) and (\ref{eqn:cexp}), if $f$,~$g\in\elltwo$ then (compare
\cite[(2.10)]{Par88})
\begin{align}
\langle\evec{f},Q_\ptn\evec{g}\rangle&=%
\prod_{n=0}^\infty\Bigl(1+\langle f_\ptn(n),g_\ptn(n)\rangle\Bigr)\nonumber\\
&=\exp\biggl(\sum_{n=0}^\infty\log\Bigl(1+\int_{\ptn_n}^{\ptn_{n+1}}%
\langle f(t),P_\ptn g(t)\rangle\std t\Bigr)\biggr)\nonumber\\
&\sim\exp\biggl(\sum_{n=0}^\infty\int_{\ptn_n}^{\ptn_{n+1}}\langle f(t),%
P_\ptn g(t)\rangle\std t\biggr)\label{eqn:asymp}\\
&\to\exp\Bigl(\int_0^\infty\langle f(t),g(t)\rangle\std t\Bigr)\qquad
\mbox{as }|\ptn|\to0,\nonumber
\end{align}
by Lemma~\ref{lem:cexp}, so $Q_\ptn\to\id_\fock$ weakly on $\evecs$, the
linear span of the set of exponential vectors; the asymptotic identity
(\ref{eqn:asymp}) holds because $\log(1+z)=z+O(z^2)$ as $z\to0$. Since each
$Q_\ptn$ is an orthogonal projection, strong convergence on $\evecs$, so on
$\fock$, follows.
\end{proof}

\section{Modified QS integrals}\label{sec:modint}

\begin{myquotation}
\textit{Natura abhorret vacuum.}

{\small\ -- Fran\c{c}ois Rabelais, \textit{Gargantua et Pantagruel}, Bk.~1,
Ch.~5 (1534).}
\end{myquotation}

To examine the behaviour of the discrete approximations which will be
constructed in the following sections, it is useful first to introduce a
slight extension of the iterated QS integral (QS being, of course, an
abbreviation for quantum stochastic).

\begin{notation}
Let $\ini$ be a fixed complex Hilbert space (the \emph{initial space}) and let
$\ffock:=\ini\otimes\fock$, $\bbebe:=\ini\otimes\bebe$ and
$\eevecs:=\ini\algten\evecs$. As is customary, the tensor sign will be omitted
before exponential vectors: $u\evec{f}:=u\otimes\evec{f}$.
\end{notation}

\begin{definition}
Given a Hilbert space $\hilb$, an \emph{$\hilb$-process} $X=(X_t)_{t\in\R_+}$
is a weakly measurable family of linear operators with common domain
$\hilb\algten\evecs$, \textit{i.e.},
\begin{equation}
X_t\in\lop{\hilb\algten\evecs}{\hilb\otimes\fock}\qquad\forall\,t\in\R_+
\end{equation}
and $t\mapsto\langle u\evec{f},X_t v\evec{g}\rangle$ is measurable for all
$u$,~$v\in\hilb$ and $f$,~$g\in\elltwo$.

An $\hilb$-process $X$ is \emph{vacuum adapted} if
\begin{equation}
\langle u\evec{f},X_t v\evec{g}\rangle=%
\langle u\evec{\indf{[0,t[}f},X_t v\evec{\indf{[0,t[}g}\rangle
\end{equation}
for all $t\in\R_+$, $u$,~$v\in\hilb$ and $f$,~$g\in\elltwo$. Equivalently, the
identity $X_t=(\id_\hilb\otimes\expn_t)X_t(\id_\hilb\otimes\expn_t)$ holds for
all $t\in\R_+$, where $\expn_t\in\bop{\fock}{}$ is the second quantisation of
the multiplication operator $f\mapsto\indf{[0,t[}f$ on $\elltwo$.

An $\hilb$-process $X$ is \emph{semi-vacuum-adapted} if
$(\id_\hilb\otimes\expn_t)X_t=X_t$ for all $t\in\R_+$; clearly every
vacuum-adapted process is semi-vacuum-adapted.

If $M\in\bopp\bigl(\elltwo\bigr)$ then an $\ini\otimes\mmul$-process $X$ is
\emph{$M$-integrable} if
\begin{equation}
\|X\nnabla^M\theta\|_{L^2({[0,t[};\ini\otimes\mmul\otimes\fock)}^2=%
\int_0^t\|X_s\nnabla^M_s\theta\|^2\std s<\infty%
\qquad\forall\,\theta\in\eevecs,\ t\in\R_+,
\end{equation}
where the \emph{modified gradient}
$\nnabla^M:\eevecs\to\ffock\oplus(\ini\otimes\elltwo\otimes\fock)$ is the
linear operator such that
\begin{equation}
u\evec{f}\mapsto[u\otimes\widehat{M f}]\evec{f}=%
\begin{pmatrix}u\evec{f}\\ [u\otimes M f]\evec{f}\end{pmatrix}
\end{equation}
and the definition $\nnabla^M_s u\evec{f}:=[u\otimes\widehat{M f(s)}]\evec{f}$
is extended by linearity.
\end{definition}

\begin{notation}
Let $\Delta\in\bop{\ini\otimes\mmul\otimes\fock}{}$ be the orthogonal
projection onto $\ini\otimes\mul\otimes\fock$ and
$\Delta^\perp:=\id_{\ini\otimes\mmul\otimes\fock}-\Delta$ the projection onto
its complement, $\ffock$.
\end{notation}

\begin{theorem}\label{thm:vacest}
Let $M\in\bopp\bigl(\elltwo\bigr)$. If $X$ is an $M$-integrable,
semi-vacuum-adapted $\ini\otimes\mmul$-process, there exists a unique
semi-vacuum-adapted $\ini$-process $\vint(X;M)$, the \emph{modified QS
integral of $X$}, such that, for all $t\in\R_+$,
\begin{equation}\label{eqn:vacest}
\|\vint(X;M)_t\theta\|\le %
c_t\|X\nnabla^M\theta\|_{L^2({[0,t[};\ini\otimes\mmul\otimes\fock)}\qquad%
\forall\,\theta\in\eevecs,
\end{equation}
where $c_t:=\sqrt{2\max\{t,1\}}$, and
\begin{equation}\label{eqn:ipint}
\langle u\evec{f},\vint(X;M)_t v\evec{g}\rangle=%
\int_0^t\langle[u\otimes\widehat{f(s)}]\evec{f},X_s%
\bigl([v\otimes\widehat{M g(s)}]\evec{g}\bigr)\rangle\std s
\end{equation}
for all $u$,~$v\in\ini$ and $f$,~$g\in\elltwo$.
\end{theorem}
\begin{proof}
This follows from the behaviour of the Bochner integral and the abstract It\^o
integral: for all $t\in\R_+$ and $\theta\in\eevecs$ let
\begin{equation}\label{eqn:intdef}
\vint(X;M)_t\theta:=\int_0^t \Delta^\perp X_s\nnabla^M_s\theta\std s+%
\Ito_t(\Delta X\nnabla^M\theta),
\end{equation}
where $\Ito_t$ is the It\^o integral on ${[0,t[}$ (the adjoint of the adapted
gradient). As $s\mapsto\Delta X_s\nnabla^M_s\theta$ is an adapted vector
process, \textit{i.e.},
$\Delta X_s\nnabla^M_s\theta\in\ini\otimes\mul\otimes\fock_{[0,s[}$ for
(almost) all $s\in\R_+$, this is a good definition, and the isometric nature
of the It\^o integral \cite[Proposition~3.28]{Blt06} implies that
\begin{align}
\|\vint(X;M)_t\theta\|^2&\le %
2t\int_0^t\|\Delta^\perp X_s\nnabla^M_s\theta\|^2\std s+%
2\int_0^t\|\Delta X_s\nnabla^M_s\theta\|^2\std s\nonumber\\
&\le c_t^2\int_0^t\|X_s\nnabla^M_s\theta\|^2\std s,
\end{align}
as claimed. The identity (\ref{eqn:ipint}) follows immediately and yields
semi-vacuum-adaptedness.
\end{proof}

\begin{remark}
As may be seen from (\ref{eqn:ipint}), the modified integral preserves
semi-vacuum-adaptedness but need not preserve vacuum-adaptedness. This
identity also shows that if $A\in\bop{\ini}{}$ commutes with $X$, in the sense
that
\begin{equation}
X_t(A\otimes\id_{\mmul\otimes\fock})=(A\otimes\id_{\mmul\otimes\fock})X_t%
\qquad\forall\,t\in\R_+,
\end{equation}
then $A$ commutes with $\vint(X;M)$ in the same sense:
$\vint(X;M)_t(A\otimes\id_\fock)=(A\otimes\id_\fock)\vint(X;M)_t$ for all
$t\in\R_+$.
\end{remark}

\begin{notation}
For all $n\ge1$ and $t\in\R_+$, let
\begin{equation}
\Delta_n(t):=\{\bt:=(t_1,\ldots,t_n)\in{[0,t[}^n:t_1<\cdots<t_n\}%
\subseteq\R_+^n
\end{equation}
and, given $M\in\bop{\elltwo}{}$, define
$(\nnabla^M)^n\in\lopp\bigl(\eevecs;%
L^2(\Delta_n(t);\ini\otimes\mmul^{\otimes n}\otimes\fock)\bigr)$
such that
\begin{equation}
(\nnabla^M)^n_\bt u\evec{f}:=\bigl((\nnabla^M)^n u\evec{f}\bigr)(\bt):=%
[u\otimes\widehat{M f}^{\otimes n}(\bt)]\evec{f}
\end{equation}
for all $u\in\ini$, $f\in\elltwo$ and $\bt\in\Delta_n(t)$, where
$\widehat{g}^{\otimes n}(\bt):=%
\widehat{g(t_1)}\otimes\cdots\otimes\widehat{g(t_n)}$ for all $g\in\elltwo$
and $\bt\in\R_+^n$.
\end{notation}

\begin{theorem}\label{thm:vaciter}
Let $M\in\bopp\bigl(\elltwo\bigr)$. If $n\ge1$,
$X\in\bop{\ini\otimes\mmul^{\otimes n}}{}$ and $Y$ is a locally uniformly
bounded, semi-vacuum-adapted $\C$-process then there exists a unique
semi-vacuum-adapted $\ini$-process $\vint^n(X\otimes Y;M)$, the \emph{modified
$n$-fold QS integral}, such that, for all $t\in\R_+$,
\begin{equation}\label{eqn:vacnorm}
\|\vint^n(X\otimes Y;M)_t\theta\|^2\le c_t^{2n}%
\int_{\Delta_n(t)}\|(X\otimes Y_{t_1})(\nnabla^M)^n_\bt\theta\|^2\std\bt
\end{equation}
for all $\theta\in\eevecs$ and
\begin{multline}
\langle u\evec{f},\vint^n(X\otimes Y;M)_t v\evec{g}\rangle\\
=\int_{\Delta_n(t)}\langle u\otimes\widehat{f}^{\otimes n}(\bt),%
X[v\otimes\widehat{M g}^{\otimes n}(\bt)]\rangle%
\langle\evec{f},Y_{t_1}\evec{g}\rangle\std\bt
\end{multline}
for all $u$,~$v\in\ini$ and $f$,~$g\in\elltwo$.
\end{theorem}
\begin{proof}
If $n=1$ then the result follows by applying Theorem~\ref{thm:vacest} to the
process $X\otimes Y:t\mapsto X\otimes Y_t$. Now suppose the theorem holds for
a particular $n\ge1$, let $X\in\bop{\ini\otimes\mmul^{\otimes n+1}}{}$ and
define $X':=\widetilde{R}_{n+1}^* X\widetilde{R}_{n+1}$, where the unitary map
$\widetilde{R}_{n+1}:\ini\otimes\mmul^{\otimes n+1}\to%
\ini\otimes\mmul^{\otimes n+1}$ implements the permutation
\begin{equation*}
u\otimes x_1\otimes x_2\otimes\cdots\otimes x_{n+1}\mapsto %
u\otimes x_2\otimes\cdots\otimes x_{n+1}\otimes x_1.
\end{equation*}
By replacing $\ini$ with $\ini\otimes\mmul$, this assumption yields a
semi-vacuum-adapted $\ini\otimes\mmul$-process $\vint^n(X'\otimes Y;M)$
such that, for all $t\in\R_+$, $u$,~$v\in\ini$, $x$,~$y\in\mmul$ and
$f$,~$g\in\elltwo$,
\begin{align*}
&\langle(u\otimes x)\evec{f},%
\vint^n(X'\otimes Y;M)_t(v\otimes y)\evec{g}\rangle\\
&\quad=%
\int_{\Delta_n(t)}\langle u\otimes x\otimes\widehat{f}^{\otimes n}(\bt),%
X'[v\otimes y\otimes\widehat{M g}^{\otimes n}(\bt)]\rangle%
\langle\evec{f},Y_{t_1}\evec{g}\rangle\std\bt\\
&\quad=%
\int_{\Delta_n(t)}\langle u\otimes\widehat{f}^{\otimes n}(\bt)\otimes x,%
X[v\otimes\widehat{M g}^{\otimes n}(\bt)\otimes y]\rangle%
\langle\evec{f},Y_{t_1}\evec{g}\rangle\std\bt.
\end{align*}
Letting
\begin{equation}
\vint^{n+1}(X\otimes Y;M):=\vint(\vint^n(X'\otimes Y;M);M)
\end{equation}
gives the result: if $t\in\R_+$, $u$,~$v\in\ini$ and $f$,~$g\in\elltwo$ then
\begin{align*}
&\langle u\evec{f},\vint^{n+1}(X\otimes Y;M)_t v\evec{g}\rangle\\
&=\int_0^t\langle[u\otimes\widehat{f(s)}]\evec{f},%
\vint^n(X'\otimes Y;M)_s\bigl([v\otimes\widehat{M g(s)}]\evec{g}\bigr)%
\rangle\std s\\
&=\int_0^t\int_{\Delta_{n}(s)}\langle u\otimes%
\widehat{f}^{\otimes n+1}(\bt,s),X[v\otimes%
\widehat{M g}^{\otimes n+1}(\bt,s)]\rangle\langle\evec{f},%
Y_{t_1}\evec{g}\rangle\std\bt\std s\\
&=\int_{\Delta_{n+1}(t)}%
\langle u\otimes\widehat{f}^{\otimes n+1}(\bt),%
X[v\otimes\widehat{M g}^{\otimes n+1}(\bt)]\rangle\langle\evec{f},%
Y_{t_1}\evec{g}\rangle\std\bt;
\end{align*}
the norm estimate (\ref{eqn:vacnorm}) (and $M$-integrability of
$\vint^n(X'\otimes Y;M)$) may be shown similarly.
\end{proof}

\begin{proposition}
If $n\ge1$, $X\in\bop{\ini\otimes\mmul^{\otimes n}}{}$ and $X\otimes\expn$ is
the vacuum-adapted $\ini\otimes\mmul^{\otimes n}$-process given by setting
$(X\otimes\expn)_t:=X\otimes\expn_t$ then
\begin{equation}
\vint^n(X):=\vint^n(X\otimes\expn;\id_\elltwo)
\end{equation}
is a vacuum-adapted, bounded process: each $\vint^n(X)_t$ extends uniquely to
an element of $\bop{\ffock}{}$, the \emph{vacuum-adapted $n$-fold quantum
Wiener integral} of~$X$.
\end{proposition}
\begin{proof}
Note first that if $Z$ is a locally uniformly bounded, vacuum-adapted
$\ini\otimes\mmul$-process and $\theta\in\eevecs$ then
\begin{align}
\|Z\nnabla\theta\|^2_{L^2({[0,t[};\ini\otimes\mmul\otimes\fock)}&\le%
\|Z\|^2_{\infty,t}\int_0^t\|(\id_{\ini\otimes\mmul}\otimes\expn_s)%
\nnabla_s\theta\|^2\std s\nonumber\\
&=\|Z\|^2_{\infty,t}%
\int_0^t\bigl(\|\expn_s\theta\|^2+\|\ag_s\theta\|^2\bigr)\std s\nonumber\\
&\le\|Z\|^2_{\infty,t}(t+1)\|\theta\|^2,\label{eqn:vacnest}
\end{align}
where $\|\cdot\|_{\infty,t}$ is the essential-supremum norm on ${[0,t[}$,
$\nnabla:=\nnabla^{\id_\elltwo}$ and $\ag$ is the adapted gradient on $\ffock$
\cite[Proposition~3.27]{Blt06}. Hence $\vint(Z;\id_\elltwo)_t$ extends to a
unique element of $\bop{\ffock}{}$ for all $t\in\R_+$ and
$\vint(Z;\id_\elltwo)$ is a locally uniformly bounded, vacuum-adapted
$\ini$-process. The result now follows from the inductive construction of
$\vint^n(X)$ given in the proof of Theorem~\ref{thm:vaciter}.
\end{proof}

\begin{proposition}\label{prp:vacitb}
Let $M$,~$N\in\bopp\bigl(\elltwo\bigr)$. If $n\ge1$,
$X\in\bop{\ini\otimes\mmul^{\otimes n}}{}$, $Y$ is a locally uniformly
bounded, semi-vacuum-adapted $\C$-process and $t\in\R_+$ then
\begin{multline}
\bigl\|\bigl(\vint^n(X\otimes Y;M)_t-\vint^n(X\otimes Y;N)_t\bigr)%
u\evec{f}\bigr\|^2\\
\le2^{n-1}c_t^{2n}\|Y\|_{\infty,t}^2\|\evec{f}\|^2\sum_{m=1}^n L_m,
\end{multline}
for all $u\in\ini$ and $f\in\elltwo$, where $\|\cdot\|_{\infty,t}$ is the
essential-supremum norm on the interval $[0,t[$,
\begin{equation}
L_m:=%
\int_{\Delta_n(t)}\|X\bigl(u\otimes\widehat{M f}^{\otimes m-1}(\bt_{m)})%
\otimes[(M-N)f](t_m)\otimes\widehat{N f}^{\otimes n-m}(\bt_{(m})\bigr)\|^2%
\std\bt
\end{equation}
$\bt_{m)}:=(t_1,\ldots,t_{m-1})$ and $\bt_{(m}:=(t_{m+1},\ldots,t_n)$.
\end{proposition}
\begin{proof}
Note that, with notation as in the proof of Theorem~\ref{thm:vaciter},
\begin{multline*}
\vint^{n+1}(X\otimes Y;M)-\vint^{n+1}(X\otimes Y;N)\\
=\vint(\vint^n(X'\otimes Y;M);M)-\vint(\vint^n(X'\otimes Y;M);N)\\
+\vint\bigl(\bigl(\vint^n(X'\otimes Y;M)-%
\vint^n(X'\otimes Y;N)\bigr);N\bigr).
\end{multline*}
Now use induction, together with (\ref{eqn:vacest}), (\ref{eqn:vacnorm}) and
the fact that
\begin{multline}
\|(\vint(Z;M)_t-\vint(Z;N)_t)u\evec{f}\|^2\\
\le c_t^2\int_0^t\|Z_s\bigl([u\otimes[(M-N)f](s)]\evec{f}\bigr)\|^2\std s,
\end{multline}
by (\ref{eqn:intdef}).
\end{proof}

\begin{definition}
If $A$ is an ordered set and $n\ge1$ then $A^{n,\uparrow}$ is the collection
of strictly increasing $n$-tuples of elements of $A$. Given $\ptn\in\Ptn$, let
\begin{equation}
{[\ptn_\bp,\ptn_{\bp+1}[}:=%
\{\bt\in\R_+^n:\ptn_{p_i}\le t_i<\ptn_{p_i+1}\ (i=1,\ldots,n)\}
\end{equation}
for all $\bp=(p_1,\ldots,p_n)\in\Z_+^{n,\uparrow}$ and, for all $t\in\R_+$, let
\begin{equation}
\Delta_n^\ptn(t):=\bigcup_{\bp\in\{0,\ldots,m-1\}^{n,\uparrow}}%
{[\ptn_\bp,\ptn_{\bp+1}[}\qquad\mbox{if }t\in{[\ptn_m,\ptn_{m+1}[}.
\end{equation}
\end{definition}

\begin{theorem}\label{thm:vacsubo}
Let $M\in\bopp\bigl(\elltwo\bigr)$. If $n\ge1$, $\ptn\in\Ptn$,
$X\in\bop{\ini\otimes\mmul^{\otimes n}}{}$ and $Y$ is a locally uniformly
bounded, semi-vacuum-adapted $\C$-process then there exists a unique
semi-vacuum-adapted $\ini$-process $\vint^n(X\otimes Y;M)^\ptn$, the
\emph{modified $n$-fold QS integral subordinate to $\ptn$}, such that, for all
$t\in\R_+$,
\begin{equation}
\|\vint^n(X\otimes Y;M)^\ptn_t\theta\|^2\le c_t^{2n}%
\int_{\Delta_n^\ptn(t)}\|(X\otimes Y_{t_1})(\nnabla^M)^n_\bt\theta\|^2\std\bt
\end{equation}
for all $\theta\in\eevecs$ and
\begin{multline}
\langle u\evec{f},\vint^n(X\otimes Y;M)^\ptn_t v\evec{g}\rangle\\
=\int_{\Delta_n^\ptn(t)}\langle u\otimes\widehat{f}^{\otimes n}(\bt),%
X[v\otimes\widehat{M g}^{\otimes n}(\bt)]\rangle%
\langle\evec{f},Y_{t_1}\evec{g}\rangle\std\bt
\end{multline}
for all $u$,~$v\in\ini$ and $f$,~$g\in\elltwo$.
\end{theorem}
\begin{proof}
When $n=1$, apply Theorem~\ref{thm:vacest} to the process
$X\otimes Y:t\mapsto X\otimes Y_t$ and let
\[
\vint^1(X\otimes Y;M)^\ptn_t:=%
\sum_{m=0}^\infty\tfn{t\in{[\ptn_m,\ptn_{m+1}[}}\vint(X\otimes Y;M)_{\ptn_m}.
\]
Now suppose the theorem holds for a particular $n\ge1$, let
$X\in\bop{\ini\otimes\mmul^{\otimes n+1}}{}$ and define
$X':=\widetilde{R}_{n+1}^* X\widetilde{R}_{n+1}$ as in the proof of
Theorem~\ref{thm:vaciter}. The semi-vacuum-adapted $\ini\otimes\mmul$-process
$\vint^n(X'\otimes Y;M)^\ptn$ is such that, for all $t\in\R_+$,
$u$,~$v\in\ini$, $x$,~$y\in\mmul$ and $f$,~$g\in\elltwo$,
\begin{multline*}
\langle(u\otimes x)\evec{f},%
\vint^n(X'\otimes Y;M)^\ptn_t(v\otimes y)\evec{g}\rangle\\
=%
\int_{\Delta_n^\ptn(t)}\langle u\otimes\widehat{f}^{\otimes n}(\bt)\otimes x,%
X[v\otimes\widehat{M g}^{\otimes n}(\bt)\otimes y]\rangle%
\langle\evec{f},Y_{t_1}\evec{g}\rangle\std\bt,
\end{multline*}
so, as
\begin{equation}
\{0,\ldots,m-1\}^{n+1,\uparrow}=%
\bigcup_{k=0}^{m-1}\bigl\{(p_1,\ldots,p_n,k):%
\bp\in\{0,\ldots,k-1\}^{n,\uparrow}\bigr\}
\end{equation}
and therefore
$\Delta_{n+1}^\ptn(\ptn_m)=\bigcup_{k=0}^{m-1}%
\bigl(\Delta_n^\ptn(\ptn_k)\times{[\ptn_k,\ptn_{k+1}[}\bigr)$,
letting
\begin{equation}
\vint^{n+1}(X\otimes Y;M)^\ptn:=%
\sum_{m=0}^\infty\tfn{t\in[\ptn_m,\ptn_{m+1}[}%
\vint(\vint^n(X'\otimes Y;M)^\ptn;M)_{\ptn_m}
\end{equation}
gives the result.
\end{proof}

\begin{proposition}\label{prp:suboiter}
Let $M\in\bop{\elltwo}{}$. If $n\ge1$, $\ptn\in\Ptn$,
$X\in\bop{\ini\otimes\mmul^{\otimes n}}{}$ and $Y$ is a locally uniformly
bounded, semi-vacuum-adapted $\C$-process then
\begin{multline}
\bigl\|\bigl(\vint^n(X\otimes Y;M)_t-\vint^n(X\otimes Y;M)^\ptn_t\bigr)%
\theta\bigr\|^2\\
\le 2^{n-1}c_t^{2n}\int_{\Delta_n(t)\setminus\Delta_n^\ptn(t)}%
\|(X\otimes Y_{t_1})(\nnabla^M)^n_\bt\theta\|^2\std\bt
\end{multline}
for all $t\in\R_+$ and $\theta\in\eevecs$.
\end{proposition}
\begin{proof}
This follows by induction and the fact that if $t\in{[\ptn_m,\ptn_{m+1}[}$
then the set
\[
\{(\bt,s):s\in{[0,\ptn_m[},\ \bt\in\Delta_n(s)\setminus\Delta_n^\ptn(s)\}%
\cup\{(\bt,s):s\in{[\ptn_m,t[},\ \bt\in\Delta_n(s)\}
\]
is contained in $\Delta_{n+1}(t)\setminus\Delta_{n+1}^\ptn(t)$.
\end{proof}

\section{The toy integral}\label{sec:toyint}

\begin{myquotation}
I see salvation in discrete individuals

{\small\ -- Anton Chekhov, Letter to I.I.~Orlov (22nd February, 1899).}
\end{myquotation}

\begin{definition}\label{def:partint}
For all $n\in\Z_+$ let $\eehm_n:\bop{\ini\otimes\mmul}{}\to\bop{\bbebe}{}$ be
the normal $*$-homomorphism such that $B\otimes C\mapsto
B\otimes\id_{\bebe_{n)}}\otimes C\otimes P^\vac_{[n+1}$, where
\begin{equation}
P^\vac_{[n+1}:\bebe_{[n+1}\to\bebe_{[n+1};\ %
\bigotimes_{m=n+1}^\infty x_n\mapsto%
\bigotimes_{m=n+1}^\infty\langle\vac_{(m)},x_n\rangle\vac_{(m)}
\end{equation}
is the orthogonal projection onto the one-dimensional subspace of
$\bebe_{[n+1}$ spanned by the vector $\otimes_{m=n+1}^\infty\vac_{(m)}$.
\end{definition}

\begin{notation}
For all $\ptn\in\Ptn$, let
\begin{equation}
D_\ptn:=\id_\ini\otimes J_\ptn^*\Pi_\ptn:\ffock\to\bbebe;\ %
u\evec{f}\mapsto u\otimes\bigotimes_{n=0}^\infty\widehat{f_\ptn(n)}
\end{equation}
and note that, as $|\ptn|\to0$,
\begin{equation}
D_\ptn^* D_\ptn=\id_\ini\otimes\Pi_\ptn^* J_\ptn J_\ptn^*\Pi_\ptn=%
\id_\ini\otimes Q_\ptn\to\id_\ffock
\end{equation}
in the strong operator topology and
\begin{equation}\label{eqn:dcoisom}
D_\ptn D_\ptn^*=\id_\ini\otimes J_\ptn^*\Pi_\ptn\Pi_\ptn^* J_\ptn=%
\id_\bbebe,
\end{equation}
since $\Pi_\ptn$ is an isometric isomorphism and $J_\ptn$ an isometry.
\end{notation}

\begin{remark}
Let $X\in\bop{\ini\otimes\mmul}{}$ and $t\in\R_+$ be fixed. For all
$\ptn\in\Ptn$, let $n=n(\ptn)\in\Z_+$ be such that $t\in[\ptn_n,\ptn_{n+1}[$
and note that
\begin{multline}\label{eqn:bn}
\langle u\evec{f},D_\ptn^*\eehm_n(X) D_\ptn v\evec{g}\rangle\\
=\prod_{m=0}^{n-1}\bigl(1+\langle f_\ptn(m),g_\ptn(m)\rangle\bigr)%
\langle u\otimes\widehat{f_\ptn(n)},X[v\otimes\widehat{g_\ptn(n)}]\rangle
\end{multline}
for all $u$,~$v\in\ini$ and $f$,~$g\in\elltwo$. As $|\ptn|\to0$,
$\ptn_n\nearrow t$ and
\begin{align}
\prod_{m=0}^{n-1}\bigl(1+\langle f_\ptn(m),g_\ptn(m)\rangle\bigr)&=%
\langle J_\ptn^*\Pi_\ptn\evec{f},J_\ptn^*\Pi_\ptn%
\evec{\indf{[0,\ptn_n[}g}\rangle\nonumber\\
&=\langle Q_\ptn\evec{f},\expn_{\ptn_n}\evec{g}\rangle\to%
\langle\evec{f},\expn_t\evec{g}\rangle.
\end{align}
To analyse the second term in the right-hand side of (\ref{eqn:bn}), let
$X=\left(\begin{smallmatrix}E&F\\G&H\end{smallmatrix}\right)$, where
$E\in\bop{\ini}{}$, $F\in\bop{\ini\otimes\mul}{\ini}$,
$G\in\bop{\ini}{\ini\otimes\mul}$ and $H\in\bop{\ini\otimes\mul}{}$. Then
\begin{align}
\langle u\otimes\widehat{f_\ptn(n)},X[v\otimes\widehat{g_\ptn(n)}]\rangle&=%
\langle u,E v\rangle+\langle u,F[v\otimes g_\ptn(n)]\rangle+%
\langle u\otimes f_\ptn(n),G v\rangle\nonumber\\
&\quad+\langle u\otimes f_\ptn(n),H[v\otimes g_\ptn(n)]\rangle;%
\label{eqn:xcomp}
\end{align}
this equation shows the necessity of scaling the components of $X$ in order to
obtain non-trivial limits. Replacing $X$ by $X_{\ptn,n}$, where
\begin{align}
\begin{pmatrix}E&F\\[0.5ex]G&H\end{pmatrix}_{\ptn,n}&:=%
\begin{pmatrix}(\ptn_{n+1}-\ptn_n)E&(\ptn_{n+1}-\ptn_n)^{1/2}F\\[0.5ex]
(\ptn_{n+1}-\ptn_n)^{1/2}G&H\end{pmatrix}\\
&\phantom{:}=%
\begin{pmatrix}(\ptn_{n+1}-\ptn_n)^{1/2}&0\\[0.5ex]0&1\end{pmatrix}%
\begin{pmatrix}E&F\\[0.5ex]G&H\end{pmatrix}%
\begin{pmatrix}(\ptn_{n+1}-\ptn_n)^{1/2}&0\\[0.5ex]0&1\end{pmatrix},%
\label{eqn:scale}
\end{align}
the right-hand side of (\ref{eqn:xcomp}) becomes
\begin{align}
&\int_{\ptn_n}^{\ptn_{n+1}}\bigl(\langle u,E v\rangle+\langle u,%
F[v\otimes g(t)]\rangle\nonumber\\
&\hspace{8em}+\langle u\otimes f(t),G v\rangle+%
\langle u\otimes f(t),H[v\otimes P_\ptn g(t)]\rangle\bigr)\std t\nonumber\\
&=\int_{\ptn_n}^{\ptn_{n+1}}\bigl(\langle u\otimes\widehat{f(t)},%
X[v\otimes\widehat{g(t)}]\rangle+%
\langle u\otimes f(t),H[v\otimes(P_\ptn g-g)(t)]\rangle\bigr)\std t.%
\label{eqn:ipident}
\end{align}
\end{remark}

\begin{theorem}\label{thm:diconv}
For all $X\in\bop{\ini\otimes\mmul}{}$ and $t\in\R_+$,
\begin{equation}
\di_\ptn(X)_t:=\sum_{n=0}^\infty%
\tfn{\ptn_{n+1}\in[0,t]}D_\ptn^*\eehm_n(X_{\ptn,n})D_\ptn\to%
\vint(X)_t\quad\mbox{as $|\ptn|\to0$}
\end{equation}
strongly on $\eevecs$.
\end{theorem}
\begin{proof}
It follows from (\ref{eqn:ipident}) that $\di_\ptn(X)_t$ can be written as the
sum of a semi-vacuum-adapted QS integral and an It\^o-integral remainder
term. Let
\begin{equation}\label{eqn:dexpdef}
\expn^\ptn_t:=\sum_{n=0}^\infty\tfn{t\in[\ptn_n,\ptn_{n+1}[}\expn_{\ptn_n}%
\qquad\forall\,t\in\R_+
\end{equation}
and note that $t\mapsto X\otimes Q_\ptn\expn^\ptn_t$ is vacuum-adapted. If
$t\in[\ptn_m,\ptn_{m+1}[$ then
\begin{align}
\langle u&\evec{f},\di_\ptn(X)_t v\evec{g}\rangle\nonumber\\
&=\sum_{n=0}^{m-1}\langle Q_\ptn\evec{f},\expn_{\ptn_n}\evec{g}\rangle%
\langle u\otimes\widehat{f_\ptn(n)},%
X_{\ptn,n}[v\otimes\widehat{g_\ptn(n)}]\rangle\nonumber\\
&=\int_0^{\ptn_m}\langle[u\otimes\widehat{f(s)}]\evec{f},%
(X\otimes Q_\ptn\expn^\ptn_s)\bigr([v\otimes\widehat{g(s)}]\evec{g}\bigr)%
\rangle\std s\nonumber\\
&\quad+\int_0^{\ptn_m}\langle[u\otimes f(s)]\evec{f},%
H[v\otimes(P_\ptn g-g)(s)]\otimes Q_\ptn\expn^\ptn_s\evec{g}\rangle\std s.
\end{align}
If $\Ito_s$ denotes the abstract It\^o integral on ${[0,s[}$ then this shows
that
\begin{multline}\label{eqn:itorem}
(\di_\ptn(X)_t-%
\vint(\indf{[0,\ptn_m[}X\otimes Q_\ptn\expn^\ptn;\id_\elltwo)_t)v\evec{g}\\
=\Ito_{\ptn_m}\bigl(H[v\otimes(P_\ptn g-g)(\cdot)]%
\otimes Q_\ptn\expn^\ptn_\cdot\evec{g}\bigr);
\end{multline}
as the It\^o integral is an isometry, the norm of this quantity is bounded
above by
\[
\|H\|\,\|v\|\,\|P_\ptn g-g\|_{L^2({[0,t[};\mul)}\|\evec{g}\|\to0%
\quad\mbox{as }|\ptn|\to0.
\]
Finally, since $Q_\ptn\expn^\ptn_s\to\expn_s$ strongly as $|\ptn|\to0$ for all
$s\in\R_+$, Theorem~\ref{thm:vacest} and the dominated-convergence theorem
imply that
\begin{equation}
\vint(X)_t-\vint(\indf{[0,\ptn_m[}X\otimes Q_\ptn\expn^\ptn;\id_\elltwo)_t\\
=\vint(X\otimes(\expn-\indf{[0,\ptn_m[}Q_\ptn\expn^\ptn);\id_\elltwo)_t
\end{equation}
tends to $0$ strongly on $\eevecs$ as $|\ptn|\to0$, as required.
\end{proof}

\section{Multiple integrals}\label{sec:mulint}

\begin{myquotation}
O, thou hast damnable iteration

{\small\ -- William Shakespeare, \textit{Henry IV, Part 1},
Act~I, Scene~ii (1596).}
\end{myquotation}

\begin{remark}
For all $X\in\bop{\ini\otimes\mmul^{\otimes 2}}{}$, $t\in\R_+$ and
$\ptn\in\Ptn$, let
\begin{equation}
\di^2_\ptn(X)_t:=\sum_{n=0}^\infty\sum_{m=0}^{n-1}\tfn{\ptn_{n+1}\in[0,t]}%
D_\ptn^*\eehm_{m,n}(X_{\ptn,m,n})D_\ptn,
\end{equation}
where $\eehm_{m,n}:\bop{\ini\otimes\mmul^{\otimes 2}}{}\to\bop{\bbebe}{}$ is
the normal $*$-homomorphism such that
\begin{equation}
B\otimes C_1\otimes C_2\mapsto B\otimes\id_{\bebe_{m)}}\otimes C_1\otimes %
P^\vac_{[m+1,n)}\otimes C_2\otimes P^\vac_{[n+1},
\end{equation}
with $P^\vac_{[m+1,n)}$ and $P^\vac_{[n+1}$ the orthogonal projections onto
$\bigotimes_{k=m+1}^{n-1}\C\vac_{(k)}$ and
$\bigotimes_{k=n+1}^\infty\C\vac_{(k)}$, respectively.

To find the correct scaling for $X_{\ptn,m,n}$, note that if
\begin{equation}
\Psi[\ptn]_n:=%
\begin{pmatrix}(\ptn_{n+1}-\ptn_n)^{1/2}&0\\0&\id_{\mul}\end{pmatrix}%
\in\bop{\mmul}{}\qquad\forall\,n\in\Z_+
\end{equation}
and $Y\in\bop{\ini\otimes\mmul}{}$ then
$Y_{\ptn,n}=%
(\id_\ini\otimes\Psi[\ptn]_n)Y(\id_\ini\otimes\Psi[\ptn]_n)$,
so let
\begin{equation}
X_{\ptn,m,n}:=%
(\id_\ini\otimes\Psi[\ptn]_m\otimes\Psi[\ptn]_n)X%
(\id_\ini\otimes\Psi[\ptn]_m\otimes\Psi[\ptn]_n)
\end{equation}
for all $(m,n)\in\Z_+^{2,\uparrow}$. Having examined the case of multiplicity
two, the general case is now clear.
\end{remark}

\begin{definition}
For all $n\ge1$, $X\in\bop{\ini\otimes\mmul^{\otimes n}}{}$, $t\in\R_+$ and
$\ptn\in\Ptn$, let
\begin{equation}
\di^n_\ptn(X)_t:=\sum_{\bp\in\Z_+^{n,\uparrow}}\tfn{\ptn_{p_n+1}\in[0,t]}%
D_\ptn^*\eehm_\bp(X_{\ptn,\bp})D_\ptn,
\end{equation}
where $\eehm_\bp:\bop{\ini\otimes\mmul^{\otimes n}}{}\to\bop{\bbebe}{}$ is the
normal $*$-homomorphism such that
\begin{equation}
B\otimes C_1\otimes\cdots\otimes C_n\mapsto %
B\otimes\id_{\bebe_{p_1)}}\otimes C_1\otimes %
P^\vac_{[p_1+1,p_2)}\otimes\cdots\otimes C_n\otimes P^\vac_{[p_n+1},
\end{equation}
in which $C_m$ acts on $\mmul_{(p_m)}$ for $m=1,\ldots,n$ and
$P^\vac:x\mapsto\langle x,\vac\rangle\vac$ acts on $\mmul_{(q)}$ for all
$q\ge p_1$ such that $q\not\in\{p_1,\ldots,p_n\}$, and
\begin{equation}
X_{\ptn,\bp}:=%
(\id_\ini\otimes\Psi[\ptn]_{p_1}\otimes\cdots\otimes\Psi[\ptn]_{p_n})X%
(\id_\ini\otimes\Psi[\ptn]_{p_1}\otimes\cdots\otimes\Psi[\ptn]_{p_n}).
\end{equation}
This is the discrete analogue of the vacuum-adapted $n$-fold quantum Wiener
integral of $X$.
\end{definition}

\begin{theorem}\label{thm:iterint}
If $n\ge1$, $X\in\bop{\ini\otimes\mmul^{\otimes n}}{}$ and $t\in\R_+$ then
\begin{equation}
\di^n_\ptn(X)_t=\vint^n\bigl(X\otimes Q_\ptn\expn^\ptn;P_\ptn\bigr)^\ptn_t%
\to\vint^n\bigl(X\bigr)_t
\end{equation}
strongly on $\eevecs$ as $|\ptn|\to0$.
\end{theorem}
\begin{proof}
If $\bp\in\Z_+^{n,\uparrow}$ then, with the obvious extension of notation,
\begin{multline}
\langle u\evec{f},D_\ptn^*\eehm_\bp(X_{\ptn,\bp})D_\ptn v\evec{g}\rangle=\\
\prod_{m=0}^{p_1-1}\langle\widehat{f_\ptn(m)},\widehat{g_\ptn(m)}\rangle%
\langle u\otimes\bigotimes_{k\in\bp}\widehat{f_\ptn(k)},X_{\ptn,\bp}%
[v\otimes\bigotimes_{k\in\bp}\widehat{g_\ptn(k)}]\rangle\label{eqn:ipmint1}
\end{multline}
for all $u$,~$v\in\ini$ and $f$,~$g\in\elltwo$. Furthermore, as
\begin{equation}
\Psi[\ptn]_k\widehat{f_\ptn(k)}=\frac{1}{\sqrt{\ptn_{k+1}-\ptn_k}}%
\int_{\ptn_k}^{\ptn_{k+1}}\widehat{f(t)}\std t,
\end{equation}
it follows that
\begin{multline}
\bigl\langle u\otimes\bigotimes_{k\in\bp}\widehat{f_\ptn(k)},X_{\ptn,\bp}%
[v\otimes\bigotimes_{k\in\bp}\widehat{g_\ptn(k)}]\bigr\rangle\\
=\int_{[\ptn_\bp,\ptn_{\bp+1}[}\bigl\langle u\otimes%
\widehat{f}^{\otimes n}(\bt),X[v\otimes%
\widehat{P_\ptn g}^{\otimes n}(\bt)]\bigr\rangle\std\bt,\label{eqn:ipmint2}
\end{multline}
which gives the identity. That the limit is as claimed may be established by
writing the difference $\di^n_\ptn(X)-\vint^n(X)$ as
\begin{multline*}
\vint^n(X\otimes Q_\ptn\expn^\ptn;P_\ptn)^\ptn-%
\vint^n(X\otimes Q_\ptn\expn^\ptn;P_\ptn)\\
+\vint^n(X\otimes (Q_\ptn\expn^\ptn-\expn);P_\ptn)+%
\vint^n(X\otimes\expn;P_\ptn)-\vint^n(X\otimes\expn;\id_\elltwo)
\end{multline*}
and employing Proposition~\ref{prp:suboiter}, Theorem~\ref{thm:vaciter} and
Proposition~\ref{prp:vacitb}.
\end{proof}

\section{Product formulae}\label{sec:itof}

\begin{myquotation}
\textit{Entia non sunt multiplicanda praeter necessitatem.}

{\small\ -- William of Ockham.}
\end{myquotation}

\begin{remark}
Given $\bp=(p_1,\ldots,p_m)\in\Z_+^{m,\uparrow}$ and
$\bq=(q_1,\ldots,q_n)\in\Z_+^{n,\uparrow}$ with $p_m<q_1$, let
\begin{equation}
\bp\cup\bq:=(p_1,\ldots,p_m,q_1,\ldots,q_n)\in\Z_+^{m+n,\uparrow}.
\end{equation}
If the normal $*$-homomorphism
\[
\fflip_{n,m}:%
\bop{\ini\otimes\mmul^{\otimes n}\otimes\mmul^{\otimes m}}{}\to%
\bop{\ini\otimes\mmul^{\otimes m}\otimes\mmul^{\otimes n}}{}
\]
is that determined by the transposition
$B\otimes C\otimes D\mapsto B\otimes D\otimes C$ then, letting
\begin{align}
Y\rhd X&:=[\fflip_{n,m}(Y\otimes\id_{\mmul^{\otimes m}})]%
(X\otimes(P^\vac)^{\otimes n})\\
\mbox{and}\quad X\lhd Y&:=(X\otimes(P^\vac)^{\otimes n})%
[\fflip_{n,m}(Y\otimes\id_{\mmul^{\otimes m}})],
\end{align}
it is readily verified that
\begin{equation}
\eehm_\bq(Y)\eehm_\bp(X)=\eehm_{\bp\cup\bq}(Y\rhd X)\quad\mbox{and}\quad
\eehm_\bp(X)\eehm_\bq(Y)=\eehm_{\bp\cup\bq}(X\lhd Y)
\end{equation}
for all $X\in\bop{\ini\otimes\mmul^{\otimes m}}{}$ and
$Y\in\bop{\ini\otimes\mmul^{\otimes n}}{}$. Furthermore, for any
$\ptn\in\Ptn$,
\begin{equation}
(Y\rhd X)_{\ptn,\bp\cup\bq}=Y_{\ptn,\bq}\rhd X_{\ptn,\bp}\quad\mbox{and}%
\quad(X\lhd Y)_{\ptn,\bp\cup\bq}=X_{\ptn,\bp}\lhd Y_{\ptn,\bq},
\end{equation}
since $\Psi[\ptn]_p C\Psi[\ptn]_p P^\vac=\Psi[\ptn]_p C P^\vac\Psi[\ptn]_p$
for all $p\in\Z_+$ and $C\in\bop{\mmul}{}$. The following Proposition is an
immediate consequence of these observations.
\end{remark}

\begin{proposition}[Fubini]
If $X\in\bop{\ini\otimes\mmul^{\otimes m}}{}$ and
$Y\in\bop{\ini\otimes\mmul^{\otimes n}}{}$ then
\begin{align}
\di^{m+n}_\ptn(Y\rhd X)_t&=%
\sum_{\bq\in\Z_+^{n,\uparrow}}\tfn{\ptn_{q_n+1}\in[0,t]}D_\ptn^* %
\eehm_\bq(Y_{\ptn,\bq})D_\ptn\di^m_\ptn(X)_{\ptn_{q_1}}\label{eqn:lfub}\\
\mbox{and}\quad\di^{m+n}_\ptn(X\lhd Y)_t&=%
\sum_{\bq\in\Z_+^{n,\uparrow}}\tfn{\ptn_{q_n+1}\in[0,t]}%
\di^m_\ptn(X)_{\ptn_{q_1}}D_\ptn^*\eehm_\bq(Y_{\ptn,\bq})D_\ptn%
\label{eqn:rfub}
\end{align}
for all $\ptn\in\Ptn$ and $t\in\R_+$.
\end{proposition}

\begin{theorem}[Quantum It\^o product formula]
If $X$,~$Y\in\bop{\ini\otimes\mmul}{}$ then
\begin{equation}\label{eqn:pito}
\vint^1(Y)\vint^1(X)=%
\vint^2(Y\rhd X)+\vint^2(Y\lhd X)+\vint^1(Y\Delta X),
\end{equation}
where $\Delta\in\bop{\ini\otimes\mmul}{}$ denotes the orthogonal projection
onto $\ini\otimes\mul$.
\end{theorem}
\begin{proof}
Note first that if $\alpha_{\ptn,n}:=(\ptn_{n+1}-\ptn_n)^{1/2}$ for all
$\ptn\in\Ptn$ and $n\in\Z_+$ then
\begin{equation}
(\id_\ini\otimes\Psi[\ptn]_n)^2=\id_\ini\otimes\Psi[\ptn]_n^2=%
\Delta+\alpha_{\ptn,n}^2\Delta^\perp,
\end{equation}
whence
\begin{align}
Y_{\ptn,n}X_{\ptn,n}=%
(X\Delta Y)_{\ptn,n}+\alpha_{\ptn,n}^2(X\Delta^\perp Y)_{\ptn,n}.
\end{align}
This working, the fact that $\eehm_m$ is a homomorphism for all $m\in\Z_+$ and
the identities (\ref{eqn:lfub}--\ref{eqn:rfub}) with $m=n=1$ imply that
\begin{equation}\label{eqn:prod}
\di_\ptn(Y)_t\di_\ptn(X)_t=\di^2_\ptn(Y\rhd X)_t+%
\di^2_\ptn\bigl(Y\lhd X)_t+\di_\ptn(Y\Delta X)_t+Z^\ptn_t
\end{equation}
for all $\ptn\in\Ptn$ and $t\in\R_+$, where
\begin{equation}
Z^\ptn_t:=\sum_{m=0}^\infty\tfn{\ptn_{m+1}\in[0,t]}\alpha_{\ptn,m}^2%
D_\ptn^*\eehm_m\bigl((Y\Delta^\perp X)_{\ptn,m}\bigr)D_\ptn.
\end{equation}
Working as in the proof of Theorem~\ref{thm:diconv} (compare
(\ref{eqn:itorem})) shows that
\begin{multline}
Z^\ptn_t u\evec{f}=%
\vint(Y\Delta^\perp X\otimes W^\ptn;\id_\elltwo)_{\ptn_n}u\evec{f}\\
+\Ito_{\ptn_n}\bigl(V[u\otimes(P_\ptn f-f)(\cdot)]\otimes %
W^\ptn_\cdot\evec{f}\bigr)
\end{multline}
for all $t\in{[\ptn_n,\ptn_{n+1}[}$, $u\in\ini$ and $f\in\elltwo$, where
$V:=\Delta Y\Delta^\perp X\Delta$ and
\begin{equation}
W^\ptn_t:=Q_\ptn\sum_{n=0}^\infty\tfn{t\in[\ptn_n,\ptn_{n+1}[}%
\alpha_{\ptn,n}^2\expn_{\ptn_n}\qquad\forall\,t\in\R_+.
\end{equation}
Now $W^\ptn_t\to0$ in norm as $|\ptn|\to0$, since
$\|W^\ptn_t\|\le\sup_{n\ge1}\alpha_{\ptn,n}^2$ for all $t\in\R_+$, so
$Z_{\ptn,n}\to0$ strongly on $\eevecs$, by Theorem~\ref{thm:vacest} and It\^o
isometry. Combining this with Theorem~\ref{thm:iterint}, it follows that
\begin{equation}
\di_\ptn(Y)_t\di_\ptn(X)_t\to\vint^2(Y\rhd X)_t+\vint^2(Y\lhd X)_t+%
\vint^1(Y\Delta X)_t
\end{equation}
strongly on $\eevecs$ as $|\ptn|\to0$ and this gives the result.
\end{proof}

\begin{remark}
The quantum It\^o formula (\ref{eqn:pito}) may be compared to that valid for
the usual form of adaptedness \cite[Exercise after Proposition~3.20]{Lin05}.
\end{remark}

\section{Further development}\label{sec:future}

\begin{myquotation}
Unbounded hopes were placed on each successive extension

{\small\ -- George Bernard Shaw, \textit{Socialism: Principles and Outlook},
Shavian Tract No.~4, The Illusions of Socialism and Socialism: Principles and
Outlook (1956).}
\end{myquotation}

This section contains little analysis, but sets out the basic situation once
one moves beyond bounded integrands.

\begin{definition}
An \emph{admissible triple} $(\ini_0,\mul_0,S)$ is a dense subspace
$\ini_0\subseteq\ini$, a dense subspace $\mul_0\subseteq\mul$ and a subset
$S\subseteq\elltwo$ such that
\begin{mylist}
\item[(i)] each $f\in S$ has compact support,
\item[(ii)] $f(t)\in\mul_0$ for all $t\in\R_+$ and $f\in S$
\item[and (iii)] $\evecs_S:=\lin\{\evec{f}:f\in S\}$ is dense in $\fock$.
\end{mylist}
\end{definition}

\begin{definition}
If $X\in\lop{\ini_0\algten\mmul_0}{\ini\otimes\mmul}$, where $\ini_0$ is a
subspace of $\ini$, $\mul_0$ is a subspace of $\mul$ and
$\mmul_0:=\C\oplus\mul_0\subseteq\mmul$, then
\begin{equation}
\eehm_n(X):=%
U_n^*(X\algten\id_{\bebe_{n)}}\algten P^\vac_{[n+1})U_n\in%
\lopp\Bigl(\ini_0\algten\bigalgten_{m=0}^\infty\mmul_0;\bbebe\Bigr)
\end{equation}
for all $n\in\Z_+$, where the unitary operator $U_n:\bbebe\to\bbebe$ is such
that
\begin{equation}
u\otimes\bigotimes_{m=0}^\infty x_m\mapsto u\otimes x_n%
\otimes\bigotimes_{m=0}^{n-1}x_m\otimes\bigotimes_{m=n+1}^\infty x_m
\end{equation}
and
\begin{equation}
\bigalgten_{m=0}^\infty\mmul_0:=%
\lin\Bigl\{\bigotimes_{m=0}^\infty x_m\Big|x_m\in\mmul_0\ %
\forall\,m\ge0,\ \exists\,l\in\Z_+:x_l=x_{l+1}=\cdots=\vac\Bigr\}.
\end{equation}
\end{definition}

\begin{proposition}
Let $(\ini_0,\mul_0,S)$ be admissible. If
$X\in\lop{\ini_0\algten\mmul_0}{\ini\otimes\mmul}$ and $t\in\R_+$ are such
that $\int_0^t\|X[u\otimes\widehat{f(s)}]\|^2\std s<\infty$ for all
$u\in\ini_0$ and $f\in S$ then
\begin{equation}
\di_\ptn(X)_t:=\sum_{n=0}^\infty%
\tfn{\ptn_{n+1}\in[0,t]}D_\ptn^*\eehm_n(X_{\ptn,n})D_\ptn\to%
\vint(X)_t
\end{equation}
weakly on $\eevecs_S:=\ini_0\algten\evecs_S$ as $|\ptn|\to0$, where
$\vint(X):=\vint(X\algten\expn;\id_\elltwo)$.
\end{proposition}
\begin{proof}
Note that
\begin{multline}
\langle u\evec{f},D_\ptn^*\eehm_n(X) D_\ptn v\evec{g}\rangle\\
=\prod_{m=0}^{n-1}\langle \widehat{f_\ptn(m)},\widehat{g_\ptn(m)}\rangle%
\int_{\ptn_n}^{\ptn_{n+1}}\langle u\otimes%
\widehat{P_\ptn f(s)},X[v\otimes\widehat{g(s)}]\rangle\std s
\end{multline}
for all $n\in\Z_+$, $u$,~$v\in\ini$ and $f$,~$g\in S$, so if
$t\in{[\ptn_n,\ptn_{n+1}[}$ then, as $|\ptn|\to0$,
\begin{align*}
\langle u\evec{f},\di_\ptn(X)_t v\evec{g}\rangle&=%
\int_0^{\ptn_n}\langle u\otimes\widehat{P_\ptn f(s)},%
X[v\otimes\widehat{g(s)}]\rangle\langle\evec{f},%
Q_\ptn\expn^\ptn_s\evec{g}\rangle\std s\\
&\to\int_0^t\langle u\otimes\widehat{f(s)},X[v\otimes\widehat{g(s)}]\rangle%
\langle\evec{f},\evec{\indf{[0,s[}g}\rangle\std s.
\end{align*}
\end{proof}

\begin{remark}
Similarly, if
$X\in\lop{\ini_0\algten\mmul_0^{\algten n}}{\ini\otimes\mmul^{\otimes n}}$ and
$t\in\R_+$ are such that
\begin{equation}
\int_{\Delta_n(t)}\bigl\|X%
[u\otimes\widehat{f}^{\otimes n}(\bt)]\bigr\|^2\std\bt<\infty%
\qquad\forall\,u\in\ini_0,\ f\in S
\end{equation}
then $\di^n_\ptn(X)_t\to\vint^n(X)_t$ weakly on $\eevecs_S$ as $|\ptn|\to0$;
the proper definitions of $\di_\ptn^n(X)$ and $\vint^n(X)$ should be clear
from the above.
\end{remark}


\section*{Bibliography}

\begin{thebibliography}{99}

\bibitem{AcB89}
\textsc{L.~Accardi} and \textsc{A.~Bach}, \textit{Central limits of squeezing
operators}, in: Quantum probability and applications IV (Rome, 1987),
L.~Accardi and W.~von Waldenfels (eds.), Lecture Notes in Mathematics~1396,
Springer, Berlin, 1989, 7--19.

\bibitem{Att03}
\textsc{S.~Attal}, \textit{Approximating the Fock space with the toy Fock
space}, S\'eminaire de Probabilit\'es~XXXVI, J.~Az\'ema, M.~\'Emery, M.~Ledoux
and M.~Yor (eds.), Lecture Notes in Mathematics~1801, Springer, Berlin, 2003,
477--491.

\bibitem{AtP06}
\textsc{S.~Attal} and \textsc{Y.~Pautrat}, \textit{From repeated to continuous
quantum interactions}, Ann.\ Henri Poincar\'e~7 (2006), 59--104.

\bibitem{Blt06}
\textsc{A.C.R.~Belton}, \textit{Some self-adjoint quantum semimartingales},
Proc.\ London Math.\ Soc.~(3)~92 (2006), 791--816.

\bibitem{Blt07}
\textsc{A.C.R.~Belton}, \textit{Random-walk approximation to vacuum cocycles},
%\href{http://arxiv.org/abs/math.OA/0702700}{\texttt{arXiv:math.OA/0702700}},
\url{arXiv:math.OA/0702700},
version~2, 2007.

\bibitem{BHJ06}
\textsc{L.~Bouten}, \textsc{R.~van Handel} and \textsc{M.R.~James},
\textit{A discrete invitation to quantum filtering and feedback control},
%\href{http://arxiv.org/abs/math.PR/0606118}{\texttt{arXiv:math.PR/0606118}},
\url{arXiv:math.PR/0606118},
version~4, 2006.

\bibitem{Bru02}
\textsc{T.A.~Brun}, \textit{A simple model of quantum trajectories},
Amer.\ J.\ Phys.~70 (2002), 719--737.

\bibitem{FrS07}
\textsc{U.~Franz} and \textsc{A.~Skalski}, \textit{Approximation of quantum
L\'evy processes by quantum random walks},
%\href{http://arxiv.org/abs/math.FA/0703339}{\texttt{arXiv:math.FA/0703339}},
\url{arXiv:math.FA/0703339},
version~1, 2007.

\bibitem{Gou04}
\textsc{J.~Gough}, \textit{Holevo-ordering and the continuous-time limit for
open Floquet dynamics}, Lett.\ Math.\ Phys.~67 (2004), 207--221.

\bibitem{GoS04}
\textsc{J.~Gough} and \textsc{A.~Sobolev}, \textit{Stochastic Schr\"odinger
equations as limit of discrete filtering}, Open Sys.\ Inf.\ Dyn.~11 (2004),
235--255.

\bibitem{KaR97}
\textsc{R.V.~Kadison} and \textsc{J.R.~Ringrose}, \textit{Fundamentals of the
theory of operator algebras. Volume II: Advanced theory}, Graduate Studies in
Mathematics 16, American Mathematical Society, Providence, Rhode Island, 1997.

\bibitem{Kum06}
\textsc{B.~K\"ummerer}, \textit{Quantum Markov processes and applications in
physics}, in: Quantum independent increment processes II, M.~Sch\"urmann and
U.~Franz (eds.), Lecture Notes in Mathematics~1866, Springer, Berlin, 2006,
259--330.

\bibitem{Lei01}
\textsc{M.~Leitz-Martini}, \textit{Quantum stochastic calculus using
infinitesimals}, doctoral thesis, Eberhard Karls Universit\"at T\"ubingen,
2001. \url{http://w210.ub.uni-tuebingen.de/dbt/volltexte/2002/458/}

\bibitem{Lin05}
\textsc{J.M.~Lindsay}, \textit{Quantum stochastic analysis --- an
introduction}, in: Quantum independent increment processes I, M.~Sch\"urmann
and U.~Franz (eds.), Lecture Notes in Mathematics~1865, Springer, Berlin,
2005, 181--271.

\bibitem{LiP88}
\textsc{J.M.~Lindsay} and \textsc{K.R.~Parthasarathy}, \textit{The passage
from random walk to diffusion in quantum probability. II}, Sankhy\=a Ser.~A~50
(1988), 151--170.

\bibitem{Mey86}
\textsc{P.-A.~Meyer}, \textit{\'El\'ements de probabilit\'es quantiques.
I--V}, in: S\'eminaire de Probabilit\'es XX, J.~Az\'ema and M.~Yor (eds.),
Lecture Notes in Mathematics~1204, Springer, Berlin, 1986, 186--312.

\bibitem{Mey89}
\textsc{P.-A.~Meyer}, \textit{\'El\'ements de probabilit\'es quantiques.
X. Approximation de l'oscillateur harmonique (d'apr\`es L.~Accardi et
A.~Bach)}, in: S\'eminaire de Probabilit\'es XXIII, J.~Az\'ema, P.-A.~Meyer
and M.~Yor (eds.), Lecture Notes in Mathematics~1372, Springer, Berlin, 1989,
175--182.

\bibitem{Par88}
\textsc{K.R.~Parthasarathy}, \textit{The passage from random walk to diffusion
in quantum probability}, J.\ Appl.\ Probab.~25A (1988), 151--166.

\bibitem{Pau05}
\textsc{Y.~Pautrat}, \textit{From Pauli matrices to quantum It\^o formula},
Math.\ Phys.\ Anal.\ Geom.~8 (2005), 121--155.

\bibitem{Sah05}
\textsc{L.~Sahu}, \textit{Quantum random walks and their convergence},
%\href{http://arxiv.org/abs/math.OA/0505438}{\texttt{arXiv:math.OA/0505438}},
\url{arXiv:math.OA/0505438},
version~1, 2005.

\bibitem{Sin06}
\textsc{K.B.~Sinha}, \textit{Quantum random walk revisited}, in: Quantum
Probability, M.~Bo\.zejko, W.~M{\l}otkowski and J.~Wysocza\'nski (eds.),
Banach Center Publications~73, Polish Academy of Sciences, Warsaw, 2006,
377--390.

\end{thebibliography}
\end{document}